\documentclass[11pt]{amsart}

\usepackage{amssymb}
\usepackage{amsmath}
\usepackage{amsthm}
\usepackage{stmaryrd}
\usepackage{amsfonts,mathrsfs}
\usepackage{fullpage}

\title[One metric extension]{A universal metric for the canonical bundle of a holomorphic family of projective algebraic manifolds}
\author{Dror Varolin}
\thanks{Partially supported by an NSF grant}
\address{
Department of Mathematics \newline \indent
Stony Brook University \newline \indent
Stony Brook, NY 11794}

\newcommand{\noi}{\noindent}

\newcommand{\co}{{\mathcal O}}

\newcommand{\si}{{\mathscr I}}

\newcommand{\sx}{{\mathscr X}}

\newcommand{\vp}{\varphi} 
\newcommand{\ve}{\varepsilon}

\newcommand{\C}{{\mathbb C}}
\newcommand{\D}{{\mathbb D}}

\newcommand{\N}{{\mathbb N}}

\newcommand{\R}{{\mathbb R}}

\newcommand{\di}{\partial}
\newcommand{\dbar}{\bar \partial}

\newcommand{\ii}{\sqrt{-1}}

\newcommand{\tensor}{\otimes}

\begin{document}
\maketitle

\theoremstyle{plain}
\newtheorem{thm}{\sc Theorem}
\newtheorem{lem}[thm]{\sc Lemma}
\newtheorem{d-thm}[thm]{\sc Theorem}
\newtheorem{prop}[thm]{\sc Proposition}
\newtheorem{cor}[thm]{\sc Corollary}

\theoremstyle{definition}
\newtheorem{conj}[thm]{\sc Conjecture}
\newtheorem{defn}[thm]{\sc Definition}
\newtheorem{qn}[thm]{\sc Question}

\theoremstyle{definition}
\newtheorem{ex}[thm]{\sc Example}

\theoremstyle{remark}
\newtheorem*{rmk}{\sc Remark}
\newtheorem*{rmks}{\sc Remarks}
\newtheorem*{ack}{\sc Acknowledgment}
\newtheorem*{sdefn}{\sc Definition}

\begin{center}{\small 
{\it Dedicated to M. Salah Baouendi on the occasion of his 60th birthday.}}
\end{center}

\section{Introduction}

In his celebrated work \cite{siu-1, siu-2}, Siu proved that the plurigenera of any algebraic manifold are invariant in families.  More precisely, let $\pi : \sx \to \D$ be a holomorphic submersion (i.e., $d\pi$ is nowhere zero) from a complex manifold $\sx$ to the unit disk $\D$, and assume that every fiber $\sx _t := \pi ^{-1} (t)$ is a compact projective manifold.  Then for every $m \in \N$, the function $P_m : \D \to \N$ defined by $P_m(t) := h^0(\sx _t , m K_{\sx _t})$ is constant.

Siu's approach to the problem begins with the observation that the function $P_m$ is upper semi-continuous.  Thus in order to prove that $P_m$ is continuous (hence constant) it suffices to show that given a global holomorphic section $s$ of $mK_{\sx _0}$, there is a family of global holomorphic sections $s_t$ of $\sx _t$, for all $t$ in a neighborhood of $0$, that varies holomorphically with $t$ and satisfies $s_0=s$.

To prove such an extension theorem, Siu establishes a generalization of the Ohsawa-Takegoshi Extension Theorem to the setting of complex submanifolds of a Kahler  manifold having codimension 1 and cut out by a single, bounded holomorphic function.  This theorem, which we will discuss below, requires the existence of a singular Hermitian metric on the ambient manifold having non-negative curvature current, with respect to which the section to be extended is $L^2$.  Thus in the presence of the extension theorem, the approach reduces to construction of such a metric.

The case where the fibers $\sx _t$ of our holomorphic family are of general type was treated in \cite{siu-1}.  In this setting, Siu produced a single singular Hermitian metric $e^{-\kappa}$ for $K_X$ so that every $m$-canonical section is $L^2$ with respect to $e^{-(m-1)\kappa}$.  

However, in the case where the fibers $\sx _t$ of our holomorphic family are assumed only to be algebraic, and not necessarily of general type, Siu's proof in \cite{siu-2} does not construct a single metric as in the case of general type.  Instead, Siu constructs for every section $s$ of $mK_{\sx _0}$ a singular Hermitian metric for $mK_{\sx}$ of non-negative curvature so that $s$ is $L^2$ with respect to this metric.  

\begin{sdefn}
Let $\sx \to \Delta$ be a holomorphic family of complex manifolds and $\sx _0$ the cental fiber of $\sx$.  A universal canonical metric for the pair $(\sx , \sx_0)$ is a singular Hermitian metric $e^{-\kappa}$ for the canonical bundle $K_{\sx}$ of $\sx$ such that for every global holomorphic section $s \in H^0(\sx _0 , mK_{\sx _0})$, 
\[
\int _{\sx _0} |s|^2 e^{-(m-1)\kappa} < +\infty.
\]
\end{sdefn}

The goal of this paper is to prove that for any holomorphic family $\sx \to \Delta$ of compact complex algebraic manifolds with central fiber $\sx _0$, the pair $(\sx , \sx_0)$ has a universal canonical metric having non-negative curvature current.  To this end, our main theorem is the following result.

\begin{thm}\label{main}
Let $X$ be a complex manifold admitting a positive line bundle $A \to X$, and $Z \subset X$ a smooth compact complex submanifold of codimension 1.  Assume there is a subvariety $V \subset X$ not containing $Z$ such that $X-V$ is a Stein manifold.  Let $T \in H^0(X,Z)$ be a holomorphic section of the line bundle associated to $Z$, thought of as a divisor.  Let $E \to X$ be a holomorphic line bundle and denote by $K_X$ the canonical bundle of $X$.  Assume we are given singular metrics $e^{-\vp _E}$ for $E$ and $e^{-\vp _Z}$ for the line bundle associated to $Z$.

Suppose in addition that the above data satisfy the following assumptions.
\begin{enumerate}
\item[(R)]  The metrics $e^{-\vp _E}$ and $e^{-\vp _Z}$ restrict to singular metrics on $Z$.

\item[(B)]  
\[
\sup _X |T|^2e^{-\vp _Z} < +\infty.
\]

\item[(G)]  The line bundles $p(K_X+Z+E)+A$, $0 \le p \le m-1$, are globally generated, in the sense that a finite number of sections of $H^0(X,p(K_X+Z+E)+A)$ generate the sheaf $\co _X(p(K_X+Z+E)+A)$.

\item[(P)]  $\ii \di \dbar \vp _E \ge 0$ and there exists a constant $\mu$ such that $\mu \ii \di \dbar \vp _E \ge \ii \di \dbar \vp _Z$.

\item[(T)]  The singular metric $e^{-(\vp _Z+\vp _E)}|Z$ has trivial multiplier ideal:  
\[
\si (Z, e^{-(\vp _Z+\vp _E)}|Z) = \co _Z.
\]
\end{enumerate}
Then there is a metric $e^{-\kappa}$ for $K_X+Z+E$ with the following properties:
\begin{enumerate}
\item[(C)]  $\ii \di \dbar \kappa \ge 0$.

\item[(L)] For every $m > 0$ and every section $s\in H^0(Z,m(K_Z+E|Z))$, $|s|^2e^{-((m-1)\kappa +\vp _E + \vp _Z)}$ is locally integrable.

\item[(I)] For every integer $m > 0$ and every section $s \in H^0(Z,m(K_Z+E))$, 
\[
\int _Z |s|^2e^{-(m-1)\kappa + \vp _E} < +\infty.
\]
\end{enumerate}
\end{thm}

\begin{rmks} 
\begin{itemize}
\item[(i)]  For the ambient manifold $X$, we have in mind the following two examples:  either $X$ is compact complex projective (in which case the variety $V$ could be taken to be a hyperplane section of some embedding of $X$) or else $X$ is a family of compact complex algebraic manifolds.  In the former case, it is well-known that the hypothesis (G) holds for any sufficiently ample $A$, while in the latter case, one might have to shrink $X$ a little to obtain (G).  Of course, there are many other examples of such $X$. 

\item[(ii)]
Note that in condition (L), the local functions $|s|^2e^{-((m-1)\kappa +\vp _E + \vp _Z)}$ depend on the local trivializations of the line bundles in question.  However, the local integrability condition is independent of these choices.
\end{itemize}
\end{rmks}

Together with a variant of the Ohsawa-Takegoshi Theorem (Theorem \ref{ot-thm} below), Theorem \ref{main} implies a generalization of Siu's extension theorem to the case where the normal bundle of the submanifold $Z$ is not necessarily trivial.  The first extension theorem of this type was established by Takayama \cite[Theorem 4.1]{tak}under some additional hypotheses.  The general case was done in \cite{v-tak}, where Theorem \ref{ot-thm} was also established.  The argument here is related to that of \cite{v-tak}, but the focus is on construction of the metric rather than on the extension theorem.

As a result of Theorem \ref{main}, we have the following corollary, which is our stated goal.

\begin{cor}\label{pseudo-thm}
For every holomorphic family $\sx \to \Delta$ of smooth projective varieties with central fiber $\sx _0$, the pair $(\sx , \sx _0)$ has, perhaps after slightly shrinking the family, a universal canonical metric having non-negative curvature current.
\end{cor}

\begin{proof}
Let $X$ be a family of compact projective manifolds $\pi : \sx \to \D$, and $Z= \sx _0$ the central fiber. Take $T=\pi$, $E = \co _{\sx}$ and $\vp _E \equiv 0$.  Since $\sx _0$ is cut out by a single holomorphic function, the line bundle associated to $\sx _0$ is trivial.  Take $\vp _Z \equiv 0$.  Then the hypotheses of Theorem \ref{main} are satisfied, perhaps after shrinking the family, and we obtain a metric $e^{-\kappa}$ for $K_{\sx}$ such that $\ii \di \dbar \kappa  \ge 0$ and $|s|^2 e^{-(m-1) \kappa _m}$ is integrable for every integer $m >0$ and every section $s \in H^0(\sx _0, mK_{\sx _0})$. 
\end{proof}

\begin{rmk}
Note that in the setting of families, the constant $\mu$ is not needed, and the hypotheses (L) and (I) are the same.
\end{rmk}

\begin{rmk}
In his paper \cite{ts}, Tsuji has claimed the existence of a metric with the properties stated in Corollary \ref{pseudo-thm}.  As in our approach, Tsuji's proof makes use of an infinite process.  It seems that convergence of this process was not checked; in fact, it is demonstrated in \cite{siu-2} that Tsuji's process, as well as any reasonable modification of it, diverges.
\end{rmk}

\begin{prop}\label{extension-metric}
For each integer $m > 0$, fix a basis $s_1 ^{(m)}, ..., s_{N_m} ^{(m)}$ of $H^0(X,m(K_Z+E|Z))$.  Choose constants $\ve _m$ such that the metric 
\[
\kappa _0 :=\log \left ( \sum _{m=1} ^{\infty} \ve _m \left ( \sum _{\ell=1} ^{N_m} |s_{\ell}^{(m)}|^2 \right )^{1/m} \right )
\]
is convergent.  Suppose $e^{-\vp _E}$ is locally integrable.  Then for each $m>0$ and every $s \in H^0(X,m(K_Z+E|Z))$, 
\[
\int _Z |s|^2 e^{-((m-1)\kappa_0 +\vp _E)} < +\infty.
\]
\end{prop} 

\begin{proof}
Fix $s \in H^0(X,m(K_Z+E|Z))$, and let $\kappa _{0,m} = \log \left ( \sum _{\ell =1} ^{N_m} |s_{\ell}^{(m)}|^2 \right )^{1/m}$.  Note that $e^{-\kappa _0} \lesssim e^{-\kappa _{0,m}}$, and thus we have 
\begin{eqnarray*}
\int _Z |s|^2 e^{-(m-1)\kappa_0 + \vp _E} &\lesssim &\int _Z |s|^2 e^{-(m-1)\kappa_{0,m} + \vp _E}\\
 &=& \int _Z |s|^{2/m} \left ( \frac{|s|^{2}}{|s^{(m)}_1|^2 + ... +|s^{(m)}_{N_m}|^2}\right )^{(m-1)/m} e^{\gamma _E- \vp _E} e^{-\gamma _E}\\
&\lesssim & \int _Z |s|^{2/m} e^{\gamma _E - \vp _E} e^{-\gamma _E}\\
&\lesssim & \left (  \int _Z |s|^2 e^{\gamma _E - \vp _E} e^{- m\gamma _E}\omega ^{-(n-1)(m-1)} \right )^{1/m} \left ( \int _Z e^{\gamma _E - \vp _E} \omega ^{n-1} \right )^{(m-1)/m},
\end{eqnarray*}
where $\omega$ is a fixed K\"ahler form for $Z$ and $e^{-\gamma _Z}$ is a smooth metric for $E|Z$.  The last inequality is a consequence of H\"older's Inequality.  Since $e^{-\vp_E}$ is locally integrable, we are done. 
\end{proof}

A calculation similar to the proof of Proposition \ref{extension-metric} shows that $|s|^2e^{-((m-1)\kappa_0 + \vp _Z + \vp _E)}$ is locally integrable on $Z$.  Thus in view of Proposition \ref{extension-metric}, Theorem \ref{main} follows if we construct a metric $e^{-\kappa}$ with non-negative curvature current such that $e^{-\kappa}|Z= e^{-\kappa _0}$.  This is precisely what we do.  We employ a technical simplification, due to Paun \cite{p}, of Siu's original idea of extending metrics using an Ohsawa-Takegoshi-type extension theorem for sections.

\tableofcontents

\section{The Ohsawa-Takegoshi Extension theorem}

Let $Y$ be a K\"ahler manifold of complex dimension $n$.  Assume there exists an analytic hypersurface $V \subset Y$ such that $Y-V$ is Stein.  Examples of such manifolds are Stein manifolds (where $V$ is empty) and projective algebraic manifolds (where one can take $V$ to be the intersection of $Y$ with a projective hyperplane in some projective space in which $Y$ is embedded).

Fix a smooth hypersurface $Z \subset Y$ such that $Z \not \subset V$.  
In \cite{v-tak} we proved the following generalization of the Ohsawa-Takogoshi Extension Theorem.

\begin{thm}\label{ot-thm}
Suppose given a holomorphic line bundle $H \to Y$ with a singular Hermitian metric $e^{-\psi}$, and a singular Hermitian metric $e^{-\vp _Z}$ for the line bundle associated to the divisor $Z$, such that the following properties hold. 
\begin{enumerate}
\item[(i)] The restrictions $e^{-\psi}|Z$ and $e^{-\vp_Z}|Z$ are singular metrics. 
\item[(ii)] There is a global holomorphic section $T \in H^0(Y,Z)$ such that 
\[
Z= \{ T= 0\} \quad \text{and} \quad \sup_Y |T|^2e^{-\vp _Z} =1.
\]
\item[(iii)] $\ii \di \dbar \psi \ge 0$ and there is an integer $\mu > 0$ such that $\mu \ii \di \dbar \psi \ge \ii \di \dbar \vp _Z$.
\end{enumerate}
Then for every $s \in H^0 (Z, K_Z + H)$ such that 
\[
\int _Z |s|^2 e^{-\psi} < +\infty \quad \text{and} \quad s\wedge dT \in \si (e^{-(\vp _Z +\psi)}|Z),
\]
there exists a section $S \in H^0 (Y,K_Y +Z + H)$ such that 
\[
S|Z = s \wedge dT \quad \text{and} \quad \int _Y |S|^2 e^{-(\vp _Z + \psi)} \le 40 \pi  \mu \int _Z |s|^2 e^{-\psi}.
\]
\end{thm}

\section{Inductive construction of certain sections by extension}

Fix a holomorphic line bundle $A \to X$ such that the property (G) in Theorem \ref{main} holds.

\noi Let us fix bases 
\[
\{\tilde \sigma ^{(m,0,p)}_j \ ;\ 1 \le j \le M_p \}
\]
of $H^0(X, p(K_{X}+Z+E)+ A )$.  We let $\sigma ^{(m,0,p)} _j \in H^0(Z,p(K_{Z} + E|Z)+ A|Z)$ be such that 
\[
\tilde \sigma ^{(m,0,p)} _j |Z = \sigma ^{(m,0,p)} _j \wedge (dT) ^{\tensor p}.
\]
We also fix smooth metrics  
\[
e^{-\gamma _Z}\ \text{and}\  e^{-\gamma _E} \ \text{ for } Z \to X,\ \text{and } E\to X 
\]
respectively.  Finally, let us fix bases
\[
s_1^{(m)}, ..., s_{N_m} ^{(m)} \text{ for } H^0(X,m(K_Z+E|Z)), \quad m=1,2,...,
\]
orthonormal with respect to the singular metric $(\omega ^{-(n-1)}e^{-\gamma _E}) ^{m-1} e^{-\vp _E}$ for $(m-1)K_Z + mE|Z$.  (Since $e^{-\vp _E}$ is locally integrable, every holomorphic section is integrable with respect to this metric.)

\begin{prop}\label{paun-prop}
For each $m =1,2,...$ there exist a constant $C_m < +\infty$ and sections 
\[
\tilde \sigma ^{(m,k,p)}_{j,\ell} \in H^0(X,(km+p) (K_X+Z+E)+A)
\]
where $p=1,2,...,m-1$, $1 \le j \le M_p$, $1 \le \ell \le N_m$ and $k=1,2,...$, with the following properties.  
\begin{enumerate}
\item[(a)] $\tilde \sigma _{j,\ell} ^{(m,k,p)} | Z = (s_{\ell} ^{(m)})^{\tensor k} \tensor \sigma _j ^{(m,0,p)} \wedge (dT) ^{(km+p)}$

\item[(b)]  If $k \ge 1$, 
\[
\int _X \frac{\sum _{j=1} ^{M_0} |\tilde \sigma _{j,\ell} ^{(m,k,0)}|^2 e^{-(\gamma _Z + \gamma _E)}}{\sum _{j=1} ^{M_{m-1}} |\tilde \sigma _{j,\ell} ^{(m,k-1,m-1)}|^2}  \le C_m.
\]

\item[(c)]  For $1 \le p \le m-1$,
\[
\int _X \frac{\sum _{j=1} ^{M_p} |\tilde \sigma _{j,\ell} ^{(m,k,p)}|^2e^{-(\gamma _Z + \gamma _E)}}{\sum _{j=1} ^{M_{p-1}} |\tilde \sigma _{j,\ell} ^{(m,k,p-1)}|^2}  \le C_m.
\]

\end{enumerate}
\end{prop}

\begin{proof} (Double induction on $k$ and $p$.)  
Fix a constant $\widehat C_m$ such that the 
\[
\sup _{X} \frac{\sum _{j=1}^{M_{0}} |\tilde \sigma _{j} ^{(m,0,0)}|^2\omega ^{n(m-1)}e^{(m-1)(\gamma _Z+\gamma _E)}}{\sum _{j=1}^{M_{m-1}} |\tilde \sigma_{j} ^{(m,0,m-1)}|^2} \le \widehat C_m
\]
and 
\[
\sup _{Z} \frac{\sum _{j=1}^{M_{0}} |\sigma_{j}^{(m,0,0)}|^2\omega ^{(n-1)(m-1)}e^{(m-1)\gamma _E}}{\sum _{j=1}^{M_{m-1}} |\sigma_{j}^{(m,0,m-1)}|^2} \le \widehat C_m, 
\]
and  for all $0 \le p \le m-2$, 
\[
\sup _X  \frac{\sum _{j=1}^{N_{p+1}} |\tilde \sigma_{j}^{(m,0,p+1)}|^2\omega ^{-n}e^{-(\gamma _Z+\gamma _E)}}{\sum _{j=1}^{M_p} |\tilde \sigma_{j}^{(m,0,p)}|^2} \le \widehat C_m, 
\]
and
\[
\sup _Z  \frac{\sum _{j=1}^{N_{p+1}} |\sigma_{j}^{(m,0,p+1)}|^2\omega ^{-(n-1)}e^{-\gamma _E}}{\sum _{j=1}^{M_p} |\sigma_{j}^{(m,0,p)}|^2} \le \widehat C_m.
\]

\noi ($k=0$) We set $\tilde \sigma ^{(m,0,p)}_{j,\ell} := \tilde \sigma ^{(m,0,p)}_j$  and simply observe that
\[
\int _X \frac{\sum _{j=1} ^{M_p} |\tilde \sigma_{j,\ell}^{(m,0,p)}|^2e^{-(\gamma _Z + \gamma _E)}}{\sum _{j=1} ^{M_{p-1}} |\tilde \sigma_{j,\ell}^{(m,0,p-1)}|^2}  \le \widehat C_m \int _X \omega ^n.
\]

\noi ($k \ge 1$)  Assume the result has been proved for $k-1$.

\noi \underline{(($p=0$))}:   Consider the sections $(s_{\ell}^{(m)})^{\tensor k} \tensor \sigma _j ^{(m,0,0)}$, and define the semi-positively curved metric 
\[
\psi _{k,\ell,0} := \log \sum _{j=1} ^{M_{m-1}} |\tilde \sigma_{j,\ell}^{(m,k-1,m-1)}|^2 
\]
for the line bundle $(mk-1)(K_X+Z+E)+A$.  Observe that locally on $Z$, 
\begin{eqnarray*}
|(s^{(m)}_{\ell}\wedge dT ^m)^k\tensor \sigma _j ^{(m,0,0)}|^2 e^{-(\vp _Z+\psi_{k,\ell,0}+\vp _E)}&=& |s^{(m)}_{\ell} \wedge dT ^m|^2 \frac{|\sigma_{j} ^{(m,0,0)}|^2 e^{-(\vp_Z + \vp _E)}}{\sum _{j=1} ^{M_{m-1}} |\sigma_{j}^{(m,0,m-1)}|^2} \\
&\lesssim& |s^{(m)}_{\ell}|^2e^{-(\vp_Z + \vp _E)}.
\end{eqnarray*}
Moreover, we have 
\[
\ii \di \dbar (\psi _{k,\ell,0}+\vp _E) \ge 0 \quad \text{and}\quad \mu \ii \di \dbar (\psi _{k,\ell,0} + \vp _E) \ge \ii \di \dbar \vp _Z .
\]
Finally,
\begin{eqnarray*}
&& \int _{Z} |(s^{(m)}_{\ell})^k\tensor \sigma _{j} ^{(m,0,0)}|^2 e^{-(\psi_{k,\ell,0}+\vp _E)} \\
&=& \int _{Z} |s^{(m)}_{\ell}|^2 \frac{|\sigma_{j}^{(m,0,0)}|^2 e^{(m-1) \gamma _E} e^{-((m-1) \gamma _E + \vp _E)}}{\sum _{j=1} ^{M_{m-1}} |\sigma_{j}^{(m,0,m-1)}|^2 }< +\infty.
\end{eqnarray*}
We may thus apply Theorem \ref{ot-thm} to obtain sections 
\[
\tilde \sigma_{j,\ell}^{(m,k,0)}\in H^0(X, mk (K_X +Z+E)+A), \quad 1 \le j \le M_0, \ 1\le \ell \le N_m,
\]
such that 
\[
\tilde \sigma_{j,\ell}^{(m,k,0)} |Z = (s^{(m)}_{\ell})^{\tensor k} \tensor \sigma_{j,\ell}^{(m,0,0)}\wedge (dT)^{\tensor km} , \quad 1 \le j \le M_0,\ 1\le \ell \le N_m,
\]
and 
\[
\int _{X} |\tilde \sigma_{j,\ell}^{(m,k,0)}|^2 e^{-(\psi_{k,\ell,0}+\vp_Z + \vp_E)} \le 40\pi \mu \int _{Z} |s^{(m)}_{\ell}|^2 \frac{|\sigma ^{(0)} _j|^2 e^{-(\vp _E + \vp _B)}}{\sum _{j=1} ^{N_{m-1}} |\sigma ^{(m-1)}_j|^2 }.
\]
Summing over $j$, we obtain
\begin{eqnarray*}
&& \int _{X} \frac{\sum _{j=1} ^{M_0} |\tilde \sigma_{j,\ell}^{(m,k,0)}|^2 e^{-(\gamma_Z + \gamma _E)}} {\sum _{j=1} ^{M_{m-1}} |\tilde \sigma_{j,\ell}^{(m,k-1,m-1)}|^2} \\
&\le& \sup _X e^{\vp _Z + \vp _E - \gamma_Z - \gamma _E} \int _{X} \frac{\sum _{j=1} ^{M_0} |\tilde \sigma_{j,\ell}^{(m,k,0)}|^2 e^{-(\vp _Z + \vp _E)}} {\sum _{j=1} ^{M_{m-1}} |\tilde \sigma_{j,\ell}^{(m,k-1,m-1)}|^2}\\
&\le&40 \pi \sup _X e^{\vp _Z + \vp _E - \gamma_Z - \gamma _E} \int _{Z} |s^{(m)}_{\ell}|^2 \frac{\sum _{j=1} ^{M_0} |\sigma_{j}^{(m,0,0)}|^2e^{-\vp _E}}{\sum _{j=1} ^{M_{m-1}}|\sigma_{j}^{(m,0,m-1)}|^2 } e^{-\kappa}\\
&\le&40\pi \widehat C _m \sup _X e^{\vp _Z + \vp _E - \gamma_Z - \gamma _E} \int _{Z} |s^{(m)}_{\ell}|^2 \omega ^{-(n-1)(m-1)}e^{-((m-1) \gamma _E + \vp _E)}\\
&=& 40 \pi \widehat C _m\sup _X e^{\vp _Z + \vp _E - \gamma_Z - \gamma _E}.
\end{eqnarray*} 

\noi \underline{(($1 \le p\le m-1$))}:  Assume that we have obtained the sections $\tilde \sigma_{j,\ell}^{(m,k,p-1)}, \ 1\le j \le M_{p-1},\ 1 \le \ell \le N_m$.  Consider the non-negatively curved singular metric 
\[
\psi _{k,\ell,p} := \log \sum _{j=1} ^{M_{p-1}} |\tilde \sigma ^{(m,k,p-1)}_{j,\ell} |^2
\]
for $(km +p-1) (K_X +Z+E)+A$.  We have 
\[
|(s^{(m)}_{\ell})^k\tensor \sigma _j ^{(m,0,p)}|^2 e^{-(\vp _Z +\psi_{k,\ell,p} + \vp _E)} = \frac{|\sigma_{j}^{(m,0,p)}|^2e^{-(\vp _Z+\vp _E)}}{\sum _{j=1} ^{M_{p-1}} |\sigma_{j}^{(m,0,p-1)}|^2 } \lesssim e^{-(\vp _Z+\vp _E)},
\]
which is locally integrable on $Z$ by the hypothesis (T).  Next, 
\begin{eqnarray*}
\int _{Z} |(s_{\ell} ^{(m)})^k\tensor \sigma_{j}^{(m,0,p)}|^2 e^{-(\psi_{k,\ell,p} + \vp _E)} &=&  \int _{Z} \frac{|\sigma_{j}^{(m,0,p)}|^2e^{-\vp _E}}{\sum _{j=1} ^{M_{p-1}} |\sigma_{j}^{(m,0,p-1)}|^2 }\\
&\le & C^{\star} \int _Z  e^{\gamma _Z} \frac{|\sigma_{j}^{(m,0,p)}|^2e^{-(\vp _Z+\vp _E)}}{\sum _{j=1} ^{M_{p-1}} |\sigma_{j}^{(m,0,p-1)}|^2}< +\infty,
\end{eqnarray*}
where 
\[
C^{\star} := \sup _Z e^{\vp _Z - \gamma _Z}.
\]
Moreover, 
\[
\ii \di \dbar (\psi _{k,\ell,p} +\vp _E) \ge 0 \quad \text{and} \quad \ii \di \dbar (\psi _{k,\ell,p} +\vp _E) \ge\ii \di \dbar \vp _Z.
\]
By Theorem \ref{ot-thm} there exist sections 
\[
\tilde \sigma_{j,\ell}^{(m,k,p)} \in H^0(X, (mk+p)(K_X +Z+E)+A), \quad 1 \le j \le M_0
\]
such that 
\[
\tilde \sigma_{j,\ell}^{(m,k,p)} |Z = (s_{\ell} ^{(m)})^{\tensor k} \tensor \sigma_{j,\ell}^{(m,0,p)} \wedge (dT)^{\tensor km +p} , \quad 1 \le j \le M_p,
\]
and 
\[
\int _{X} |\tilde \sigma_{j,\ell}^{(m,k,p)}|^2 e^{-(\psi_{k,\ell,p}+\vp _Z +\vp _E)} \le 40 \pi \mu \int _{Z}\frac{|\sigma_{j}^{(m,0,p)}|^2e^{-\vp _E}}{\sum _{j=1}^{M_{p-1}} |\sigma_{j}^{(m,0,p-1)} |^2 }.
\]
Summing over $j$, we obtain
\[
\int _{X} \frac{\sum _{j=1} ^{M_p} |\tilde \sigma_{j,\ell}^{(m,k,p)}|^2e^{-(\gamma _Z +\gamma _E)}} {\sum _{j=1} ^{M_{p-1}} |\tilde \sigma_{j,\ell}^{(m,k,p-1)}|^2} \le 40\pi \mu \sup _X e^{\vp _Z + \vp _E - \gamma _Z - \gamma _E} \widehat C_m \int _{Z}e^{-\vp _E} \omega ^{n-1}.
\]
Letting 
\[
C_m := 40\pi \mu \widehat C _m \max \left ( \int _X \omega ^n,\ \sup _X e^{\vp _Z + \vp _E + \vp _B - \gamma_Z - \gamma _E},\ \sup _X e^{\vp _Z + \vp _E - \gamma _Z - \gamma _E}  \int _{Z}e^{-\vp _E} \omega ^{n-1}\right )
\]
completes the proof.
\end{proof}

\section{Construction of the metric}

\subsection{A metric associated to $\mathbf{m(K_X+Z+E)}$}

Fix a smooth metric $e^{-\psi}$ for $A \to X$.  Consider the functions 
\[
\lambda ^{(m)}_{\ell,N} := \log \sum _{j=1} ^{M_p} |\tilde \sigma_{j,\ell}^{(m,k,p)}|^2\omega ^{-n(mk+p)} e^{-(km (\gamma _Z + \gamma_E) + \psi)},
\]
where $N=mk+p$.  Set
\[
\lambda ^{(m)} _{N} := \log \sum _{\ell =1} ^{N_m} e^{\lambda ^{(m)}_{\ell, N}}.
\]
\begin{lem}\label{jensen}
For any non-empty open subset $V \subset X$ and any smooth function $f: \overline{V} \to \R _+$,
\begin{eqnarray*}
\frac{1}{\int _V f \omega ^{n}} \int _V (\lambda^{(m)} _N - \lambda ^{(m)} _{N-1}) f \omega ^{n} \le \log \left ( \frac{N_m C_m\sup _V f}{\int _V f \omega ^{n}} \right ).
\end{eqnarray*}
\end{lem}

\begin{proof}
Observe that by Proposition \ref{paun-prop}, there exists a constant $C_m$ such that for any open subset $V \subset X$,
\[
\int _V (e^{\lambda ^{(m)}_{\ell,N} - \lambda ^{(m)}_{\ell,N-1}}) f \omega ^{n} \le C_m \sup _V f, 
\]
and thus 
\[
\int _V (e^{\lambda ^{(m)}_{N} - \lambda ^{(m)}_{N-1}}) f \omega ^{n} = \sum _{\ell =1} ^{N_m} \int _V (e^{\lambda ^{(m)}_{\ell,N} - \lambda ^{(m)}_{\ell,N-1}}) f \omega ^{n} \le N_mC_m \sup _V f.
\]

An application of (the concave version of) Jensen's inequality to the concave function $\log$ then gives 
\[
\frac{1}{\int _V f\omega ^n } \int _V (\lambda ^{(m)}_N - \lambda ^{(m)} _{N-1})f \omega ^n \le \log \left ( \frac{N_mC_m \sup _V f}{\int _V f \omega ^n }\right ).
\]
The proof is complete.
\end{proof}

Consider the function 
\[
\Lambda ^{(m)} _k = \frac{1}{k} \lambda^{(m)} _{mk}.
\]
Note that $\Lambda ^{(m)}_k$ is locally the sum of a plurisubharmonic function and a smooth function.  By applying Lemma \ref{jensen} and using the telescoping property, we see that for any open set $V \subset X$ and any smooth function $f : \overline{V} \to \R _+$, 
\begin{equation}\label{L-bound}
\frac{1}{\int _V f \omega ^n } \int _V \Lambda ^{(m)} _k f \omega ^n\le m \log \left ( \frac{N_m C_m\sup _V f}{\int _V f\omega ^n } \right ) .
\end{equation}

\begin{prop}\label{lambda-prop}
There exists a constant $C^{(m)}_o$ such that 
\[
\Lambda^{(m)} _k (x) \le C^{(m)}_o, \quad x \in X.
\]
\end{prop}

\begin{proof}
Let us cover $X$ by coordinate charts $V_1,...,V_N$ such that for each $j$ there is a biholomorphic map $F_j$ from $V_j$ to the ball $B(0,2)$ of radius $2$ centered at the origin in $\C ^n$, and such that if $U_j = F^{-1}_j(B(0,1))$, then $U_1,...,U_N$ is also an open cover.  Let $W_j = V_j \setminus F^{-1}_j  (B(0,3/2))$.

Now, on each $V_j$, $\Lambda ^{(m)} _k$ is the sum of a plurisubharmonic function and a smooth function.  Say $\Lambda ^{(m)} _k = h + g$ on $V_j$, where $h$ is plurisubharmonic and $g$ is smooth.  Then for constant $A_j$ we have 
\begin{eqnarray*}
\sup _{U_j} \Lambda^{(m)}_k&\le&\sup _{U_j} g  + \sup _{U_j} h \\
&\le& \sup _{U_j} g  + A_j  \int_{W_j}h \cdot  F_{j*} dV\\
&\le& \sup _{U_j} g - A_j \int_{W_j}g \cdot  F_{j*} dV + A_j \int_{W_j}\Lambda^{(m)}_k \cdot  F_{j*} dV\\
\end{eqnarray*}
Let 
\[
C^{(m)}_j := \sup _{U_j} g - A_j  \int_{W_j}g \cdot  F_{j*} dV
\]
and define the smooth function $f_j$ by 
\[
f_j \omega ^n = F_{j*} dV.
\]
Then by \eqref{L-bound} applied with $V= W_j$ and $f=f_j$, we have
\[
\sup _{U_j} \Lambda^{(m)} _k \le C^{(m)}_j + mA_j \log \left ( \frac{N_mC_m\sup _{W_j} f_j}{\int _{W_j} f_j \omega ^n} \right ) \int _{W_j} f_j \omega ^n.
\]
Letting 
\[
C^{(m)}_o := \max _{1\le j \le N} \left \{ C^{(m)}_j + m A_j \log \left ( \frac{N_m C_m\sup _{W_j} f_j}{\int _{W_j} f_j \omega ^n} \right ) \int _{W_j} f_j \omega ^n \right \}
\]
completes the proof.
\end{proof}

Since the upper regularization of the lim sup of a uniformly bounded sequence of plurisubharmonic functions is plurisubharmonic (see, e.g., \cite[Theorem 1.6.2]{h}), we essentially have the following corollary.

\begin{cor}
The function 
\[
\Lambda ^{(m)} (x) := \limsup _{y \to x} \limsup _{k \to \infty} \Lambda ^{(m)}_k(y)
\]
is locally the sum of a plurisubharmonic function and a smooth function.  
\end{cor}

\begin{proof}
One need only observe that the function $\Lambda _k$ is obtained from a singular metric on the line bundle $m(K_X+Z+E)$  (this singular metric $e^{-\kappa^{(m)} _k}$ will be described shortly)  by multiplying by a fixed smooth metric of the dual line bundle. 
\end{proof}

Consider the singular Hermitian metric $e^{-\kappa^{(m)}}$ for $m(K_X+Z+E)$ defined by
\[
e^{-\kappa^{(m)}} = e^{-\Lambda^{(m)}} \omega ^{-nm}e^{-m(\gamma _Z+\gamma _E)}.
\]
This singular metric is given by the formula
\[
e^{-\kappa ^{(m)}(x)} = \exp \left (- \limsup _{y\to x} \limsup _{k \to \infty} \kappa^{(m)}_k (y) \right ),
\]
where
\[
e^{-\kappa^{(m)}_k} = e^{-\Lambda^{(m)}_k} \omega ^{-nm}e^{-m(\gamma_Z+\gamma _E)}.
\]
The curvature of $e^{-\kappa^{(m)}_k}$ is thus
\begin{eqnarray*}
\ii \di \dbar \kappa^{(m)} _k &=& \frac{\ii}{k} \di \dbar \log \sum _{\ell =1}^{N_m}\sum _{j=1} ^{N_0} |\tilde \sigma_{j,\ell} ^{(m,k,0)}|^2 - \frac{1}{k} \ii \di \dbar \psi \\
&\ge &  - \frac{1}{k} \ii \di \dbar \psi 
\end{eqnarray*}

We claim next that the curvature of $e^{-\kappa}$ is non-negative.  To see this, it suffices to work locally.  Then we have that the functions 
\[
\kappa ^{(m)}_k + \frac{1}{k} \psi
\]
are plurisubharmonic.  But
\[
 \limsup _{y\to x} \limsup _{k \to \infty}\kappa ^{(m)}_k + \frac{1}{k} \psi =  \limsup _{y\to x} \limsup _{k \to \infty} \kappa^{(m)}_k = \kappa ^{(m)}.
 \]
It follows that $\kappa^{(m)}$ is plurisubharmonic, as desired.

\subsection{The metric for $\mathbf{K_X+Z+E}$; Proof of Theorem \ref{main}}

Let $\ve _m$ be constants, chosen so $\ve _m \searrow 0$ sufficiently rapidly that the sum 
\[
e^{\kappa} := \sum _{m=1} ^{\infty} \ve _m e^{\tfrac{1}{m}\kappa ^{(m)}} = \sum _{m=1} ^{\infty} \exp ( \tfrac{1}{m} \kappa ^{(m)} + \log \ve _m).
\]
converges everywhere on $X$ (to a metric for $-(K_X+Z+E)$).  It is possible to find such constants since, by Proposition \ref{lambda-prop}, each $\kappa ^{(m)}$ is locally uniformly bounded from above.  (The lower bound $e^{\kappa ^{(m)}} \ge 0$ is trivial.)  Moreover, by elementary properties of plurisubharmonic functions, $\kappa$ is plurisubharmonic.  Indeed, for any $r\in \N$, the function 
\[
\psi _r := \log \sum _{m=1} ^{r} \exp ( \tfrac{1}{m} \kappa ^{(m)} + \log \ve _m)
\]
is plurisubharmonic, and $\psi_r \nearrow \kappa$.  It follows that $\kappa  = \sup _r \psi _r $ is plurisubharmonic.  (Again, see \cite[Theorem 1.6.2]{h}.)Thus $e^{-\kappa}$ is a singular Hermitian metric for $K_X+Z+E$ with non-negative curvature current.

Observe that, after identifying $K_Z$ with $(K_X+Z)|Z$ by dividing by $dT$, 
\[
\kappa ^{(m)} _k |Z = \log \left ( \sum _{\ell=1} ^{N_m} |s_{\ell}^{(m)}|^2  \right ) + \frac{1}{k} \log \sum _{j=1} ^{M_0} |\sigma_{j}^{(m,0,0)}|^2.
\]
Thus we obtain $e^{-\kappa^{(m)}} |Z = \left ( \sum _{\ell=1} ^{N_m} |s_{\ell} ^{(m)}|^2\right )^{-1}$.  It follows that 
\[
e^{-\kappa}|Z = \frac{1}{\sum _{m=1} ^{\infty} \ve _m \left (\sum _{\ell =1} ^{N_m} |s_{\ell}^{(m)}|^2 \right ) ^{2/m}}.
\]
In view of the short discussion following the proof of Proposition \ref{extension-metric}, the metric $e^{-\kappa}$ satisfies the conclusions of Theorem \ref{main}.  The proof of Theorem \ref{main} is thus complete.\qed

\begin{ack}
I am indebted to Lawrence Ein and Mihnea Popa.  It is to a discussion with them that the present paper owes its existence.  
\end{ack}

\end{document}